\documentclass[a4paper,12pt]{article}
\usepackage{amsmath,amsxtra,amssymb,latexsym, amscd,amsthm}
\usepackage{indentfirst}
\usepackage[mathscr]{eucal}

\newtheorem{Theorem}{Theorem}
\newtheorem{Lemma}{Lemma}

\newtheorem{Proposition}{Proposition}

\usepackage[utf8]{vietnam}
\usepackage{hyperref}

\begin{document}
\title{The uniform controllability property of semidiscrete approximations for the parabolic distributed parameter systems in Banach spaces}

% The thanks line in the title should be filled in if there is
% any support acknowledgement for the overall work to be included
% This \thanks is also used for the received by date info, but
% authors are not expected to provide this.

\author{Thuy NGUYEN\thanks{Université d'Orléans, Laboratoire Mathématiques et Applications, Physique Mathématiques d'Orléans (MAPMO), Bâtiment de Mathématiques, B.P. 6759, 45067 Orléans cedex 2, FRANCE. Email address : {\tt nthuyc8@yahoo.com}}}
%\date{\textrm{\small{(communicated with Emmanuel TRELAT)}}}   
\date{}

%%%%%%%%
\maketitle
%%%%%%%%
\renewcommand{\abstractname}{Abstract}
 \renewcommand{\proofname}{Proof}
\begin{abstract}
The problem we consider in this work is to minimize the $L^q$-norm $(q>2)$ of the semidiscrete controls. As shown in \cite{LT06}, under the main approximation assumptions that the discretized semigroup is uniformly analytic and that the degree of unboundedness of control operator is lower than 1/2, the uniform controllability property of semidiscrete approximations for the parabolic systems is achieved in $L^2$. In the present paper, we show that the uniform controllability property still continue to be asserted in $L^q  (q>2)$ even with the condition that the degree of unboundedness of control operator is greater than 1/2. Moreover, the minimization procedure to compute the approximation controls is provided.  An example of application is implemented for the one dimensional heat equation with Dirichlet boundary control.
\end{abstract}

% \begin{keywords}: Controllability; Partial differential equation; Discretization; Observability inequality; Duality
% \end{keywords}
%%%%%%%%%%%
\section{Introduction}  
\par Consider an infinite dimensional linear control system
\begin{equation}\dot{y}(t)=Ay(t)+Bu(t), \qquad  y(0)=y_0,\end{equation}
where the state y(t) belongs to a reflexive Banach space X, the control u(t) belongs to a reflexive Banach space U, $A:D(A)\to X$ is an operator, and B is a control operator (in general, unbounded) on $U$. Discretizing this partial differential equation by using, for instance, a finite difference or a finite element scheme, leads to a family of finite dimensional linear control systems
\begin{equation}\overset{.}{y_h}(t)=A_hy_h(t)+B_hu_h(t), \qquad  y_h(0)=y_{0h}, \end{equation}
\ where $y_h(t)\in X_h$ and $u_h(t)\in U_h$, for $0<h<h_0$ .
\par Let $y_1\in X$. If the control system (1) is controllable in time T then there exists a solution $y(.)$ of (1) associated with a control $u$ such that $y(T)=y_1$. As known, we have many available methods in order to identify the controllability. We refer to J.-L. Lions \cite{L88} for a well-known method in attainning the control of the minimal $L^2$ norm- the so called HUM (HUM stands for Hilbert uniqueness method). In this work, however, we investigate a method which we can achieve the minimization procedure in $L^q$ norm $(q>1)$. Namely, we will establish some conditions ensuring the existence and convergence of the discretized control of the minimal $L^q$ norm 
%The following question arises naturally: is it possible to find controls $u_h$, for $0<h<h_0$, converging to the control u as the mesh size h of the discretization process tends to zero, and such that the associated trajectories $y_h$, the solution of (2), converges to y(.)?
%\par The question was investigated by the authors of \cite{LT06}. They got the affirmative answer in the case the controllability of (1) is achieved in $L^2$ norm by using the Hilbert uniqueness method (HUM). It is the well-known method , introduced in \cite{L88} , which consists in minimizing a cost function, namely, the $L^2$ norm of the control. In this paper, we investigate the above question in the case where the minimization procedure of (1) is achieved in  $L^q$ norm($1<q<\infty$ ).  Our objective is to 
\begin{eqnarray}{\min}\frac{1}{q}{\int_0^T{\left\|{u_h(t)}\right\|^q}dt} \ \ (q>1) .\end{eqnarray}

\par  Necessary conditions for optimal control in finite dimensional state space were derived by Pontryagin et al \cite{PS62} (see also \cite{E}, \cite{T05}). The Maximum Principle as a set of necessary conditions for optimal control in infinite dimensional space has been studied by many authors. Since it is well known that the Maximum Principle may be false in infinite dimensional space, there are still many papers that give some conditions to ensure that the Maximum Principle remains true. It was Li and Yao \cite{LY85} who used the Eidelheit separation theorem and of the Uhl's theorem in order to extend the Maximum Principle to a large class of problems in infinite dimensional spaces when the target set is convex and the final time T is fixed. Additionally, the authors of  \cite{F87}, \cite{FF91}, \cite{LY91}, by making use of Ekeland's variational principle, give some conditions on the reachable set and on the target set in order to get an extension of Maximum Principle. Nevertheless, the problem is that when we applied the result of \cite{F87}, \cite{LY91} for the system (1) in the case the final state and final time are fixed, the finite-codimensional condition in \cite{F87}, \cite{LY91} does not satisfy for the system (1) in general. Hence we cannot use Maximum Principle in our problem.  Fortunately, thank to the Fenchel-Rockafellar duality theorem which is used in the same manner in \cite{CGL94}, \cite{GLH08}, the constrainted minimization of the function can be replaced by the unconstrainted minimization problem of corresponding conjugate function. Therefore, we will consider the above minimization procedure in the same framework with \cite{GLH08}.
\par The uniform controllability is an important area of control theory research and it has been the subject of many papers in recent years. The main goal of this article is to establish conditions ensuring a uniform controllability property of the family of discretized control systems (2) in $L^q$ and to establish a computationally feasible approximation method for identifying controllability.
\par  It is well known that controllabilty and observability are dual aspects of the same problem. We therefore will focus on the uniform observability which is shown to hold when the observability constant of the finite dimensional approximation systems does not depend on $h$. Some relevant references concerning this property has been investigated by many authors in series of articles \cite{IZ99}, \cite{LZ98}, \cite{LZ02}, \cite{NZ03}, \cite{Zua99}, \cite{Zua02}, \cite{Zua04}, \cite{Zua05}, \cite{Zua06}, \cite{BHL10a}, \cite{BHL10b} and \cite{LT06}. For finite difference schemes, a uniform observability property holds for one-dimensional heat equation \cite{LZ02}, beam equation \cite{LZ98}, Schrodinger equations \cite{Zua05}, but does not hold for 1-D wave equations \cite{IZ99}. This is due to the fact that the discrete dynamics generates high frequency spurious solutions for which the group velocity vanishes that do not exist at the continuous level. To overcome these high frequency spurious for wave equations, \cite{Zua05} showed some remedies such that Tychonoff's regularization, multigrid method, mixed finite element and filtering of high frequency, etc.
\par To our knowledge, in 1-D heat equation case, due to the fact that the dissipative effect of the 1-D heat equation acts as a filtering mechanicsm by itself and it is strong enough to exclude high frequency spurious oscillations\cite{LZ98}. However, the situation is more complex in multi-dimensional. The counter-example is shown in \cite{Zua06} for the simplest finite difference semi-discretization scheme for the heat equation in the square. 

\par In recent works in $L^2$-norm, by means of discrete Carleman inequalities, the authors in \cite{BHL10a}, \cite{BHL10b} obtain the weak uniform observability inequality for parabolic case by adding reminder terms of the form $e^{-Ch^{-2}}\left\|{\psi_h(T)}\right\|^2_{L^2(\Omega)}$ which vanishes asymptotically as $h\to 0$. Moreover, as in \cite{LT06}, the approximate controllability is derived from using semigroup arguments and introducing a vanishing term of the form $h^{\beta}\left\|{\psi_h(T)}\right\|^2_{L^2(\Omega)}$ for some $\beta >0$.
\par In fact, an efficient computing the null control for a numerical approximation scheme of the heat equation is itself a difficult problem. Accoding to the pioneering work of Carthel, Glowinski and Lions in \cite{CGL94}, the null control problem is reduced to the minimization of a dual conjugate function with respect to final condition of the adjoint state. However, as a consequence of high regularizing property of the heat kenel, this final condition does not belong to $L^2$, but a much large space that can hardly be approximated by standard techniques in numerical analysic. Recently, A. Munch and collaborators  have developed some feasible numericals such that the transmutation method, variational approach, dual and primal algorithms allow to more efficiently compute the null control (see in series \cite{CM09}, \cite{CM10}, \cite{MZ10}, \cite{PM10}).

\par The discretization framework in this paper is the same spirit as \cite{LT06}, \cite{LT00}. In \cite{LT06}, under standard assumptions on the discretization process and for an exactly null controllable parabolic system (1),  if \textbf{the degree of unboundedness of the control operator is lower than 1/2} then the semidiscrete approximation models are uniformly controllable and they also showed that for the (2), the minimizing of the cost function of discretized control  with power $q =2$ is obtained.
\par In this article, we prove the existence of the minimum of the cost function of discretized control power $q$ $(q>2)$ for type (2), in the case the operator A generates an analytic semigroup. Of course, due to regularization properties, the control system (1) is not exactly controllable in general. Hence, we focus on exact null controllability. Our main result, theorem 3.1, states that for exactly null controllable parabolic system (1) and under standard approximation assumptions, if the discretized semigroup is uniformly analytic, and if \textbf{the degree of unboundedness of the control operator B with respect to A is greater than 1/2}, then a uniform observability inequality in ($L^p$) is proved. We stress that we do not prove uniform exact null controllability property for the approximating system (2). Moreover, a minimization procedure to compute the approximation controls is provided.

\par The outline of the paper is as follows. In Section 2, we briefly review some well-known facts on controllability of linear partial differential equation in reflexive Banach spaces. In Section 3, we consider the existence and unique solution of the minimization problem in continuous case. By making use of the Fenchel-Rockafellar duality theorem, we gives a constructive way to build the control of minimal $L^q$ norm. The main result is stated in Section 4 and proved in Section 5.  An example of application and numerical simulations are provided in Section 6, for the one-dimensional heat equation with Dirichlet boundary control. An Appendix is devoted to the proof of a lemma.

\section{A short review on controllability of linear partial differential equations in reflexive Banach spaces}

\par In this section, it is convenient to first have a quick look to controllability of infinite dimensional linear control systems in reflexive Banach spaces (see more \cite{CT06}, \cite{P83}, \cite{TW07}).
\par The notation $L(E,F)$ stands for the set of linear continuous mappings from $E$ to $F$, where $E$ and $F$ are reflexive Banach spaces.
\par Let X be a reflexive Banach space. Denote by $<,>_X$ the inner product on $X$, and by $\left\|{.}\right\|_X$ the associated norm. Let $S(t)$ denote a strongly continuous semigroup on $X$, of generator $(A,D(A))$. Let $X_{-1}$ denote the completion of $X$ for norm $\left\|{x}\right\|_{-1}=\left\|{(\beta I-A)^{-1}x}\right\|$, where $\beta \in \rho(A)$ is fixed. Note that $X_{-1}$ does not depend on the specific value of $\beta \in \rho(A)$. The space $X_{-1}$ is isomorphic to $(D(A^*))'$, the dual space of $D(A^*)$ with respect to the pivot space X, and $X\subset X_{-1}$, with a continuous and dense embedding. The semigroup $S(t)$ extends to a semigroup on $X_{-1}$, still denoted $S(t)$, whose generator is an extension of the operator $A$, still denoted $A$. With these notations, $A$ is a linear operator from $X$ to $X_{-1}$.
   
\par Let $U$ be a reflexive Banach space. Denote by $<,>_U$ the inner product on $U$, and by  $\left\|{.}\right\|_U$ the associated norm.
\par A linear continuous operator $B: U \to X_{-1}$ is admissible for the semigroup $S(t)$ if every solution of 
\begin{equation} y'=Ay(t)+Bu(t),\end {equation}\\
with $y(0)=y_0 \in X$ and $u(.)\in L^q(0,+\infty;U)$, satisfies $y(t)\in X$, for every $t\ge 0$. The solution of equation (1) is understood in the mild sense, i.e,
\begin{equation} y(t)=S(t)y(0)+\int_0^T{S(t-s)Bu(s)ds}, \end{equation}
\\ for every $t\ge 0$.
\par For $T>0$, define $L_T : L^q(0,T;U)\to X_{-1} $ by    
\begin{equation} L_T u = \int_0^T {S(T-s)Bu(s)ds}. \end{equation}

\par A control operator $B\in L(U,X_{-1})$ is admissible, if and only if $Im L_T \subset X$, for some (and hence for every) $T>0$.
\par The adjoint $L^*_T $ of  $L_T$ satisfies
\begin{eqnarray}&& L^*_T : X^* \to  (L^q(0,T;U))^*=L^p(0,T;U^*) \nonumber\\
&& L^*_T \psi(t)=B^*S(T-t)^*\psi \end{eqnarray}, a.e on [0,T] for every $\psi \in D(A^*)$. Moreover, we have
\begin{equation} \left\|{L^*_T \psi}\right\| = \underset{\left\|u\right\|_q \le 1}{sup}\int_0^T {<B^*S(T-s)^*\psi,u(s)>ds}, \end{equation} for every $\psi\in X^*$ .
\par Let $B\in L(U,X_{-1})$ denote an admissible control operator.

\par We use two following lemmas (for the proofs we refer to \cite{TW07})
\begin{Lemma} $Z, X$ are reflexive Banach spaces. $G\in L(Z,X)$ then the following statements are equivalent:
\begin{itemize}
\item $G$ is onto.
\item $G^*$ bounded from below i.e there exists $C>0$ such that
\begin{equation} \left\|{G^*x}\right\|_Z \ge C\left\|{x}\right\|_X    ~\textrm{every}~  x\in X .\nonumber\end{equation}
\end{itemize}
\end{Lemma}

\begin{Lemma}$Z_1,Z_2,Z_3$ are reflexive Banach spaces. And $f\in L(Z_1,Z_3)$ , $g\in L(Z_2,Z_3)$. Then the following statements are equivalent:
\begin{itemize}
\item $Im f \subset Im g$.
\item There exists a constant $C>0$ such that :  $\left\|{f^*z}\right\|_{Z_1}\le C\left\|{g^*z}\right\|_{Z_2}$ for every $z\in Z_3$.
\item There exists an operator $h\in L(Z_1,Z_2)$ such that $f=gh$.
\end{itemize}
\end{Lemma}

\par We state some concepts as follows\\

\par For $y_0\in X$, and $T>0$, the system (4) is exactly controllable from $y_0$ in time T if for every $y_1 \in X$, there exists $u(.)\in L^q(0,T;U)$ so that the corresponding solution (4), with $y(0)=y_0$ satisfies $y(T)=y_1$.
\par In fact that the system (4) is exactly controllable  from $y_0$ in time T if and only if $L_T$ is onto, that is $Im L_T=X$. Making use of Lemma 1, there exists $C>0$ such that
\begin{eqnarray} C\left\|{\psi}\right\|_X\le \left\|{L^*_T \psi}\right\|& =& \underset{\left\|u\right\|_q \le 1}{sup}\int_0^T {<B^*S(T-s)^*\psi,u(s)>ds}\nonumber\\
& \le&  \underset{\left\|{u}\right\|_q\le 1}{sup}{\int_0^T  {\left\|{B^*S(T-t)^*{\psi}}\right\|}{\left\|{u(t)}\right\|}dt}\nonumber\\
&\le&  (\int_0^T  {\left\|{B^*S(t)^*{\psi}}\right\|^p}dt)^\frac{1}{p}.\nonumber\end{eqnarray}

\par  Therefore, the system (4) is exactly controllable in time T if and only if 
\begin{equation}  \int_0^T  {\left\|{B^*S(t)^*{\psi}}\right\|^p}dt\ge C\left\|{\psi}\right\|^p_X.\end{equation}\\

\par For $T>0$, the system (4) is said to be exactly null controllable in time T if for every $y_0\in X$, there exists $u(.)\in L^q(0,T;U)$ so that the corresponding solution of ( 4), with $y(0)=y_0$ satisfies $y(T)=0$ .
\par This means that the system (4) is exactly null controllable in time T if and only if $Im S(T)\subset Im L_T$. Making use of Lemma 2 and the same argument as above, there exists $C>0$ such that

\begin{equation} C\left\|{S(T)^*\psi}\right\|_X \le \left\|{L_T^*\psi}\right\| \le  (\int_0^T  {\left\|{B^*S(t)^*{\psi}}\right\|^p}dt)^\frac{1}{p}.\nonumber\end{equation}
\par Thus, the system (4) is exactly null controllable in time T if and only if 

\begin{equation} \int_0^T  {\left\|{B^*S(t)^*{\psi}}\right\|^p}dt \ge C\left\|{S(T)^*\psi}\right\|^p_X.\end{equation}

\section{Duality}

\par The goal of this section is to show that, with using duality arguments and Fenchel- Rockafellar theorem we can achieve the control of minimal $L^q $ norm $(q>1)$for the continuous framework.
\par Consider the system :
\begin{equation}
\begin{cases}
{\overset{.}{y}(t)}=Ay(t)+Bu(t)   \    ~\textrm{on}\ Q_T=(0,T)\times\Omega&\\
y(0)=y_0&
\end{cases}
\end{equation}\\
where $B$ is admissible and $A$ generates an analytic semigroup $S(t)$ in the reflexive Banach space $X$.
\par  Our aim is to mimimize the following functional:

\begin{equation} \begin{cases}
 ~\textrm{Minimize} \ J(u)=\frac{1}{q}\int_0^T{\left\|{u}\right\|^q dt} \ \   (q>1)&\\
 ~\textrm{Subject to} \ u\in E&
\end{cases}
\end{equation}
where $E=\left\{u\in U:{~\textrm{u steering the system from $y_0$}   ~\textrm{at time zero to y(T)=0}    }\right\}$.

\begin{Theorem} The problem (12) has a unique solution $u$.
\end{Theorem}
\par \begin {proof} 
\par First of all, we show the existence of the solution of  the optimal problem.
\par Consider a minimizing sequence $(u_n)_{n\in N}$ of controls on $[0,T]$, i.e,
\begin{equation}\int_0^T{\left\|{u_n}\right\|^q dt} ~\textrm{converges to } \inf {J(u)}~\textrm{as}\ {n\to +\infty}.
\end{equation}
\par Hence, $(u_n)_{n\in N}$ bounded in $L^q(0,T;U)$.
\par Since $U$ is reflexive space and  $q<+\infty$, then $L^q(0,T;U)$ is also reflexive.
\par Thus, up to a sequence, $u_n$ converges weakly to $u$ in $L^q$. Note that the trajectory $y_n$ (resp. $y$) associated with the control $u_n$ (resp. $u$) on [0,T] through the system
 
\begin{equation}\dot{y}_n=Ay_n+Bu_n,\ y_n(0)=y_0,\nonumber\end{equation}
and the solution of the above system is expressed in form
\begin{equation}y_n(t)=S(t)y_0+\int_0^T{S(T-s)Bu_n ds}.\nonumber
\end{equation}
\par A passage to the limit imples that
\begin{equation}\dot{y}=Ay+Bu,\ y(0)=y_0,\nonumber\end{equation}
and the solution y associated with control u in the form
\begin{equation}y(t)=S(t)y_0+\int_0^T{S(T-s)Bu ds}.\nonumber
\end{equation}
\par As $u_n$ converges weakly to $u$ in $L^q$, we get the inequality
\begin{eqnarray}\int_0^T{\left\|{u(t)}\right\|^q dt}&\le& \underset{n\to +\infty}{lim \inf}\int_0^T{\left\|{u_n(t)}\right\|^q dt}\nonumber\\
&=&\inf \int_0^T{\left\|{u(t)}\right\|^q dt}.\nonumber
\end{eqnarray}

\par It follows easily that
\begin{equation}\int_0^T{\left\|{u(t)}\right\|^q dt}=\inf \int_0^T{\left\|{u(t)}\right\|^q dt}.\nonumber
\end{equation}
\par Hence u is optimal of (12). This ensures the existence of  a optimal control. 
\par Moreover, the cost function is strictly convex then the solution is obvious uniqueness.
\end{proof}

\par By making use of convex duality, the problem of control to trajectories is reduced to the minimization of the corresponding conjugate function. Roughly speaking, it is stated through the following theorem:

\begin{Theorem}
\par (i) We have the identity:
\begin{equation}\underset{u\in E}{\inf}\frac{1}{q}\int_0^T{\left\|{u}\right\|^q dt}=-\underset{\psi_T }{\inf} (\frac{1}{p}\int_0^T{\left\|{B^*\psi}\right\|^p dt}+<\psi(0),y_0>),
\end{equation}
\ where $\psi$ be solution of :
\begin{eqnarray}&&-\dot{\psi}=A^*\psi\\
&&\psi(T)=\psi_T.
\end{eqnarray} 
\par Or, we have in the form

\begin{equation}\underset{u\in E}{\inf}\frac{1}{q}\int_0^T{\left\|{u}\right\|^q dt}=-\underset{\psi\in X^*}{\inf} (\frac{1}{p}\int_0^T{\left\|{B^*S(T-t)^*\psi}\right\|^p dt}+<S(T)^*\psi,y_0>).
\end{equation}

\par (ii) If $u_{op}$ is optimal of  the problem (12) then
\begin{equation}u_{op}(t)=\left\|{B^*S(T-t)^*\varphi}\right\|^{p-2}B^*S(T-t)^*\varphi,\nonumber
\end{equation}
\ where $\varphi$ be optimal of the function:
\begin{equation}J^*(\psi)=\frac{1}{p}\int_0^T{\left\|{B^*S(T-t)^*\psi}\right\|^p dt}+<S(T)^*\psi,y_0>. \nonumber
\end{equation}
\end{Theorem}
\par \begin {proof} \

\par (i) Let $\bar{y}$ be solution of (1) with $u=0$ and we introduce the operator $N\in L(L^q(Q_T),X)$ with $Nu=z_{u}(.,T)$ for all $u\in L^q(Q_T)$, where $z_{u}$ is solution to
\begin{eqnarray}&&\dot{z}=Az+Bu\\
&&z(x,0)=0.
\end{eqnarray}

\par Accordingly, the solution y of (11) can be decomposed in the form
\begin{equation}y=z_{u}+\bar{y}.
\end{equation}

\par The adjoint $N^*$ is given as follows
\par For each $\psi_T \in X^*$, $N^*\psi_T=B^*\psi$ where $\psi$ is solution of (15), (16).

\par Let us introduce the following functions $F$ and $G$
\begin{center}$F(z_{T})=\begin{cases}
{0} \  \ \ \ \ \ ~\textrm{for}~ z(T)=-\bar{y}(T) &\\
{+\infty} \ \ ~\textrm{otherwise}&
\end{cases}$,
\end{center}

\begin{equation}G(u)=\frac{1}{q}\int_0^T{\left\|{u}\right\|^q dt}.\nonumber
\end{equation}

\par Then, the problem (12), where the infimium is taken over all u satisfying E, is equivalent to the following minimization problem
\begin{equation}\underset{u\in L^q(Q_T)}{\inf}(F(Nu)+G(u)).
\end{equation}

\par We can apply now duality theorem of W.Fenchel and T.R.Rockafellar (see Theorem 4.2 p.60 in \cite{ET99}). It gives
\begin{equation}\underset{u\in L^q(Q_T)}{\inf}(F(Nu)+G(u))=-\underset{\psi_T\in X^*}{\inf}(G^*(N^*\psi_T)+F^*(-\psi_T)),
\end{equation} 
where $F^*$ and $G^*$ are the convex conjugate of F and G, respectively. Denote that $\psi_T=\psi(T)$, $z_T=z(T)$.
\par Note that
\begin{equation}F^*(\psi_T)=\underset{z_T=-\bar{y_T}}{\sup}<z_T,\psi_T>=-<\psi_T,\bar{y_T}>, \nonumber
\end{equation}
\ for all $\psi_T\in X^*$.

\par Additionally,
\begin{equation}G^*(\omega)=\frac{1}{p}\int_0^T{\left\|{\omega}\right\|^p dt}.\nonumber
\end{equation}

\par Therefore,
\begin{equation}G^*(N^*\psi_T)+F^*(-\psi_T)=\frac{1}{p}\int_0^T{\left\|{B^*\psi}\right\|^p dt}+<\psi_T(x),\bar{y_T}(x)>.
\end{equation} 

\par Finally, multiplying the state equation (15) by $\bar{y}$ and due to (11), we obtain
\begin{equation}
<\psi_T,\bar{y_T}>=<\psi(0),y_0>.\nonumber
\end{equation}

\par Rewrite (23) as follows
\begin{eqnarray}G^*(N^*\psi_T)+F^*(-\psi_T)&=&\frac{1}{p}\int_0^T{\left\|{B^*\psi}\right\|^p dt}+<\psi(0),y_0>\nonumber\\
&=&\frac{1}{p}\int_0^T{\left\|{B^*S(T-t)^*\psi_T}\right\|^p dt}+<S(T)^*\psi_T,y_0>,\nonumber
\end{eqnarray}
\ since $\psi$ is the solution of (15) , (16).

\par From (21) and (22), we have the identity
\begin{equation}\underset{u\in E}{\inf}\frac{1}{q}\int_0^T{\left\|{u}\right\|^q dt}=-\underset{\psi_T }{\inf} (\frac{1}{p}\int_0^T{\left\|{B^*\psi}\right\|^p dt}+<\psi(0),y_0>),\nonumber
\end{equation}
\ where $\psi$ be solution of (15), (16). 
\par Or, we have in the form

\begin{equation}\underset{u\in E}{\inf}\frac{1}{q}\int_0^T{\left\|{u}\right\|^q dt}=-\underset{\psi \in X^*}{\inf} (\frac{1}{p}\int_0^T{\left\|{B^*S(T-t)^*\psi}\right\|^p dt}+<S(T)^*\psi,y_0>).\nonumber
\end{equation}\\

\par (ii) If we denote by $(u_{op})$ , $(\varphi_{T})$ the unique solution to "LHS of (14)" and "RHS of (14)" respectively, then one has
\begin{equation}0=\frac{1}{q}\int_0^T{\left\|{u_{op}}\right\|^q dt}+\frac{1}{p}\int_0^T{\left\|{B^*\varphi_{T}}\right\|^p dt}+<\varphi_{T}(0),y_0>.
\end{equation}
\par We apply Young's inequality for the first two terms of RHS (24)
\begin{equation}\frac{1}{q}\int_{Q_T}{\left\|{u_{op}}\right\|^q dt}+\frac{1}{p}\int_{Q_T}{\left\|{B^*\varphi_{T}}\right\|^p dt}\ge \int_{Q_T}{u_{op}.B^*\varphi_{T}} .
\end{equation}

\par Then, "RHS of (24)" $\ge \int_{Q_T}{u_{op}.B^*\varphi_{T}}+<\varphi_{T}(0),y_0>$.
\par Futhermore, by multiplying two sides of (15) by y and applying Green's formula, we obtain

\begin{equation}<B^*\varphi_{T},u>+<\varphi_{T}(0),y_0>=0.
\end{equation} 

\par On the one hand, "RHS of (24)" $\ge 0$ ( due to (26)).

\par On the other hand, "RHS of (24)" $= 0$ ( due to (24)).
\par This is equivalent to that the sign "=" in inequality (14) happens, i.e,
\begin{equation}\left\|{u_{op}}\right\|^q=\left\|{B^*\varphi_{T}}\right\|^p.\nonumber
\end{equation}
\par It is also rewritten as follows
\begin{equation}u_{op}(t)=\left\|{B^*S(T-t)^*\varphi}\right\|^{p-2}B^*S(T-t)^*\varphi,\nonumber
\end{equation}
where $\varphi$ be optimal of the function $J^*$ is given as above.

\end{proof}

\par \textbf{Remark}: It is easily seen that the functional $J^*$ is convex, and from the inequality (10), is coercive. Then, it follows that  $J^*$attains a unique minimum in some point $\varphi \in D(A^*)$. As above explanation, the control ${\bar{u}}$ is chosen by
\begin{equation}\bar{u}(t)=\left\|{B^{*} S(T-t)^*{\varphi}}\right\|^{p-2}{B^{*}S(T-t)^*}{\varphi},\end{equation}
\\ for every $t\in [0,T]$ and let y(.) be the solution of (11), such that $y(0)=y_0$, associated with the control ${\bar{u}}$, then we have y(T)=0.

\par Therefore, ${\bar{u}}$ is the control of minimal of $L^q$ norm, among all controls whose associated trajectory satisfied $y(T)=0$.

\par We emphasize that observability  in $L^p$ norm ($1<p<2$) implies controllability and gives a constructive way to build the control of minimal $L^q$ norm ($q>2$). A similar result was known in $L^2$ norm through using HUM ( see more \cite{CT06}, \cite{L88}, \cite{Zua04}, \cite{Zua05}).

\section{ The main result}      
\par We are concerned in this work with the uniform controllability property for the parabolic systems. As shown in \cite{LT06}, this property is known to hold with the degree of unboundedness of control operator $\gamma\in\left[{0,1/2}\right)$. In this section, we also establish some appropriate assumptions and conditions ensuring that the unform controllability still holds in the case $\gamma\in\left[{1/2,\frac{1}{p}}\right)$.
\par Let $X$ and $U$ be Hilbert spaces, and let $  A: {D(A)} \to X$ be a linear operator and self-adjoint, generating a strongly continuous semigroup $S(t)$ on $X$.  Let $B \in L(U,D(A^*)')$ be a control operator. We make the following assumptions that will be used along this article (also refer to \cite{LT06})
\par \textbf{(H1)} The semigroup {S(t)} is analytic.
\par Therefore, (see \cite{P83}) there exist positive real number $ C_1 $ and ${ \Omega }$ such that

 \begin{equation}\left\|{S(t)}\right\|_X  \leqslant  C_1{e^{{\omega}t}}\left\|{y}\right\|_X ,  \left\|{AS(t)y}\right\|_X \leqslant   C_1\frac{{e^{{\omega}t}}}{t}\left\|{y}\right\|_X,\end{equation}
for all $t>0$ and $y\in D(A)$, and such that, if we set $\hat{A}=A-{\omega}I$, for $\theta\in \left[{0,1}\right]$  and there holds
\begin{equation}
\left\|{(-\hat{A}^{\theta})S(t)y}\right\|_X \leq C_1\frac{{e^{{\omega}t}}}{t^{\theta}}\left\|{y}\right\|_X,
\end{equation}
\\ for all $t>0$ and $y\in D(A)$.
\par Of course, inequalities (28) hold as well if one replaces $A$ by $A^{*}$, $S(t)$ by $S(t)^{*}$, for $y\in D(A^{*})$.

\par Moreover, if $\rho{(A)}$ denotes the resolvent set of A, then there exists $\delta\in \left({0,\frac{\pi}{2}}\right)$ such that
\ $\rho(A){\supset}\Delta_\delta=\left\{{\omega+\rho e^{i\theta}| \theta>0,{\left|{\theta}\right|}\leq \frac{\pi}{2}+{\delta}}\right\}$.\\

For $\lambda\in \rho (A)$, denote by $R(\lambda,A)=(\lambda I-A)^{-1}$ the resolvent of A . It follows from the previous estimates that exists $C_2>0$ such that 

\begin{equation}{\left\|{R(\lambda,A)}\right\|_{L(X)}} \leq \frac{C_2}{\left|{\lambda-\omega}\right|},   \left\|{AR(\lambda,A)}\right\|_{L(X))}\leq {C_2},\end{equation}\\
for every $\lambda\in \Delta_\delta$, and 

\begin{equation}\left\|{R(\lambda,\hat{A}))}\right\|_{L(X)}\leq \frac{C_2}{\left|{\lambda}\right|}, \left\|{\hat{A}R(\lambda,\hat{A})}\right\|_{L(X)}\leq {C_2},\end{equation}
for every $\lambda\in{\Delta_\delta+\omega}$. Similarly, inequalities (30)  and (31) hold as well with $A^{*}$ and $\hat{A}^{*}$.
\par \textbf{(H2)} The degree of unboundedness of B is $\gamma$. Assume that  $\gamma\in\left[{1/2,\frac{1}{p}}\right)$ \  ( where  p,q are conjugate, i.e $\frac{1}{p}+\frac{1}{q}=1$ and $ 1\le p<2$). This means that
\begin{equation}B\in {L(U,D((-\hat{A}^{*})^{\gamma})')}.\end{equation}
\ In these conditions, the domain of $B^{*}$ is $D(B^{*})=D((-\hat{A}^{*})^{\gamma})$, and there exists $C_3>0$ such that
\begin{equation}{\left\|{B^{*}\psi}\right\|}\leq {C_3}{\left\|{((-\hat{A}^{*})^{\gamma})\psi}\right\|}_X,\end{equation}
\ for every $\psi\in{D((-\hat{A}^{*})^{\gamma})}$.

\par \textbf{(H3)} We consider two families $(X_h)_{0<h<h_0}$ and $(U_h)_{0<h<h_0}$ of finite dimentional spaces, where h is the discretization parameter.

\par For every $h\in \left({0,h_0}\right)$, there exist the linear mappings $P_h:D((-\hat{A}^{*})^{\frac{1}{2}})'\to {X_h}$ and $\tilde{P_h}:{X_h}\to D((-\hat{A}^{*})^{\frac{1}{2}})$ and $\hat{(A^*)}^{-\gamma+\frac{1}{2}}: D(-(\hat{A^*})^{\frac{1}{2}}) \to D(-(\hat{A^*})^{\gamma})$
\  (resp., there exist linear mappings ${Q_h}:U\to {U_h} $ and $\tilde{Q_h}:{U_h}\to U$), satisfying the following requirements:

\par $(H_{3.1})$ For every $h\in \left({0,h_0}\right)$. The following properties hold

\begin{equation}{P_h}{\tilde{P_h}}=id_{X_h} ~\textrm{and}~ {Q_h}{\tilde{Q_h}}=id_{U_h}.\end{equation}

\par $\left({H_{3.2}}\right)$ There exist $s>0$ and ${C_4}>0$ such that there holds, for every $h\in \left({0,h_0}\right)$ ,
\begin{equation}\left\|{\left({I-{{\hat{(A^*)}^{-\gamma+\frac{1}{2}}}\tilde{P_h}{P_h}}}\right)\psi}\right\|_X \le {C_4}{h^s}\left\|{A^{*}\psi}\right\|_X,\end{equation}

\begin{equation}\left\|{((-\hat{A}^{*})^{\gamma})}{\left({I-{{\hat{(A^*)}^{-\gamma+\frac{1}{2}}}\tilde{P_h}{P_h}}}\right)\psi}\right\|_X \le {C_4}{h^{s\left({1-\gamma}\right)}}\left\|{A^{*}\psi}\right\|_X,\end{equation}
for every $\psi \in {D{\left({A^{*}}\right)}}$  and 
\begin{equation}\left\|{(I-\tilde{Q_h}{Q_h} )u}\right\|_U \to 0,\end{equation}
for every $u \in U$, and 
\begin{equation}\left\|{(I-\tilde{Q_h}{Q_h} )B^*{\psi}}\right\|_U \le {C_4}h^{s(1-\gamma)}{\left\|{A^*\psi}\right\|_X},\end{equation}
for every $\psi \in D(A^*)$

\par For every $h\in (0,h_0)$, the vector space $X_h$ (resp. $U_h$) is endowed with the norm $\left\|{.}\right\|_{X_h}$ (resp. $\left\|{.}\right\|_{U_h}$) defined by:
\par $\left\|{y_h}\right\|_{X_h}=\left\|{\tilde{P_h}{y_h}}\right\|_X$ for $y_h\in X_h$ (resp. ,$\left\|{u_h}\right\|_{U_h}=\left\|{\tilde{Q_h}{u_h}}\right\|_U$).
\par Therefore, we have the properties
\begin{equation}\left\|{\tilde{P_h}}\right\|_{L(X_h,X)}=\left\|{\tilde{Q_h}}\right\|_{L(U_h,U)}=1~\textrm{and}~ \left\|{{\hat{(A^*)}^{-\gamma+\frac{1}{2}}}x}\right\|_{X} \le C \left\|{x}\right\|_X,\end{equation}
\begin{equation}\left\|{P_h}\right\|_{L(X,X_h)}\le C_5~\textrm{and}~\left\|{Q_h}\right\|_{L(U,U_h)}\le C_5.\end{equation}

\par ($H_{3.3}$) For every  $h\in \left({0,h_0}\right)$, there holds
\begin{equation}{P_h}={\tilde{P_h}^{*}}~\textrm{and}~{Q_h}={\tilde{Q_h}^{*}},\end{equation}
where the adjoint operators are considered with respect to the pivot spaces $X$, $U$, $X_h$, $U_h$.

\par ($H_{3.4}$) There exists ${C_6}$ such that 
\begin{equation}\left\|{B^{*}{\hat{(A^*)}^{-\gamma+\frac{1}{2}}}\tilde{P_h}\psi_h}\right\|_U\le {C_6}h^{-\gamma s}\left\|{\psi_h}\right\|_{X_h},\end{equation}
\ for all $h\in \left({0,h_0}\right)$ and $\psi_h \in {X_h}$.
\par For every $h\in \left({0,h_0}\right)$, we define the approximation operators $A^{*}_h:{X_h}\to{X_h} $ of $A^{*}$ and $B^{*}_h:{X_h}\to{U_h} $ of $B^{*}$, by
\begin{equation}{A^{*}_h}={P_h}{A^{*}}{\tilde{P_h}}~\textrm{and}~{B^{*}_h}={Q_h}{B^{*}}{\hat{(A^*)}^{-\gamma+\frac{1}{2}}}{\tilde{P_h}}.\end{equation}

\par ($H_4$) The following properties hold:
\par ($H_{4.1}$) The family of operators $e^{tA_h}$ is uniformly analytic, in sense that there exists ${C_7}>0$ such that
\begin{equation}\left\|{e^{tA_h}}\right\|_{L(X_h)}\le {C_7}e^{\omega t},\end{equation}
\begin{equation}\left\|{ {A_h}e^{tA_h}}\right\|_{L(X_h)}\le {C_7}\frac{e^{\omega t}}{t},\end{equation}
\ for all $t>0$ and $h\in (0,h_0)$.

\par ($H_{4.2}$) There exists ${C_8}>0$ such that, for every $f\in X$ and every $h\in (0,h_0)$, the respective solutions of $\hat{A}^{*}\psi=f$ and $\hat{A}^{*}_h\psi_h={P_h}f$ satisfy
\begin{equation}\left\|{ {P_h}\psi-\psi_h}\right\|_{X_h}\le {C_8}h^s\left\|{f}\right\|_X.\end{equation}
\par In other words, there holds $\left\|{ {P_h}{\hat{A}}^{*-1}-{{\hat{A}}^{*-1}}{P_h}}\right\|_{L(X,X_h)}\le {C_8}h^s$ .\\

\par \textbf{Remark 4.1} Compare to \cite{LT06}, the important point to note here is the appearance of  the function $(\hat{A^*})^{-\gamma+\frac{1}{2}}$ in (35), (36) and (42).
\par According to \cite{LT06}, the inequality (22) make sense since $\gamma<\frac{1}{2}$ and thus $im\tilde{P_h}\subset D( (-\hat{A^*})^{1/2})\subset D( (-\hat{A^*})^{\gamma}) $.
\par In our context, on account of  $\gamma\ge \frac{1}{2}$, the inequality (36), which is similar to inequality (22) in \cite{LT06}, only make sense  if we add the functional $(\hat{A^*})^{-\gamma+\frac{1}{2}}$ in order that $im(\hat{A^*}^{-\gamma+\frac{1}{2}}\tilde{P_h})\subset D( (-\hat{A^*})^{\gamma}) $. The choice of the function  $(\hat{A^*})^{-\gamma+\frac{1}{2}}$  seems to be the best adapted to our theory.\\
\par Namely, we give here for instance about the functional $(\hat{A^*})^{-\gamma+\frac{1}{2}}$ through the heat equation with Dirichlet boundary control as follows
\begin{equation}\overset{.}{y}= \Delta y+c^2 y  ~\textrm{in}~ (0,T) \times \Omega \nonumber\end{equation}
\begin{equation} y(0,.)=y_0   ~\textrm{in}~ \Omega\nonumber\end{equation}
\begin{equation} y=u          ~\textrm{in}~  (0,T)\times \Gamma=\Sigma.\nonumber\end{equation}
\par Set $X=L^2(\Omega)$ and $U=L^2(\Gamma)$. It can be written in the form (4), where the self-adjoint operator $A:D(A)\to X$ is defined by
\begin{equation}Ay=\Delta y+c^2y:D(A)=H^2\bigcap H^1_0 \to L^2.\nonumber
\end{equation}
\par In this case, the degree of unbounded of B is $\gamma =\frac{3}{4}+\epsilon (\epsilon>0)$ (see \cite{LT00}, section 3.1).
\par We may take $\hat{A}$ as follows
\begin{equation}\hat{A}h=-\Delta h,  \ D(\hat{A})=H^2\bigcap H^1_0.\nonumber
\end{equation}
\par Therefore, $(\hat{A^*})^{-\gamma+\frac{1}{2}}=(\hat{A^*})^{-\frac{1}{4}+\epsilon}:H^1(\Omega)\to H^{\frac{3}{2}+\epsilon}(\Omega)$.\\

\par \textbf{Remark 4.2} By means of the condition of the degree of unbounded of operator B and (33), we imply that B is admissible.
\par Indeed, we have
\begin{eqnarray} \left\|{L^*_T {\psi}}\right\| &=& \underset{\left\|u\right\|_q \le 1}{sup}\int_0^T {<B^*S^*(T-s)x,u(s)>ds}\nonumber\\
&\le & (\int_0^T  {\left\|{B^*S(t)^*{\psi}}\right\|^p}dt)^\frac{1}{p}\nonumber\\
&\le&{C_3}(\int_0^T  {\left\|{{(-\hat{A^*})^{\gamma}}S(t)^*{\psi}}\right\|^p}dt)^\frac{1}{p}\nonumber\\
&\le &{C_3}(\int_0^T  {{\frac{e^{{\omega}t}}{t^{p\gamma}}}\left\|{\psi}\right\|^p}dt)^\frac{1}{p}dt \ (p\gamma<1)\nonumber\\
&\le&C_T \left\|{\psi}\right\| .\nonumber \end{eqnarray} 
 
\par \textbf{Remark 4.3}. It is easily seen that assumptions $(H_3)$ (except for the inequalities (35), (36), (42)) and $(H_{4.2})$ hold for most of the classical numerical approximation schemes, such as Galerkin methods, centered finite difference schemes,...Additionally, by using some approximation properties and the properties of the functional $(\hat{A^*})^{-\gamma+1/2}$, we prove that the inequalities (35), (36), (42) also hold for most of the above classical schemes (see the proof in Section 5). As noted in \cite{LT00}, the assumption $H_{4.1}$ of uniform analyticity is not standard, and has to be checked in each specific case.

\par The main result of the article is the following : 

\begin{Theorem}
\par Under the previous assumptions, if the control system $\dot{y}=Ay+Bu$ is exactly null controllable in time $T>0$, then there exist $\beta>0, {h_1}>0$, and positive real numbers C, C' satisfying
\begin{eqnarray}C{\left\|{e^{TA^{*}_h}\psi_h}\right\|}^p_{X_h} &\le& \int_0^T{\left\|{ {B^{*}_h}{e^{tA^{*}_h}}{\psi_h}}\right\|}^p_{U}dt+h^{\beta}{\left\|{\psi_h}\right\|}^p_{X_h}\nonumber\\
&\le & C'{\left\|{\psi_h}\right\|}^p_{X_h}, \end{eqnarray}\\
 for every $h \in ( 0,{h_1})$ and every ${\psi_h}\in {X_h}$,     $(1\le p< 2)$.

\par In these conditions, for every $y_0\in X$, and every $h\in(0,h_1)$, there exists a solution $\varphi_h\in {X_h}$ minimizing the functional
\begin{equation}  J_h(\psi_h)=\frac{1}{p}\int_0^T{\left\|{ {B^{*}_h}{e^{tA^{*}_h}}{\psi_h}}\right\|}^p_{U}dt+\frac{1}{p}{h^{\beta}}{\left\|{\psi_h}\right\|}^p_{X_h}+<e^{TA^{*}_h}\psi_h,{P_h}{y_0}>_{X_h}, (1\le p<2) \end{equation}
\ and the sequence $(\tilde{Q_h}u_h)_{0<h<h_1}$ , where the control $u_h$ is defined by
\par $$u_h(t)=\left\|{B^{*}_h e^{(T-t)A^{*}_h}{\varphi_h}}\right\|^{p-2}{B^{*}_he^{(T-t)A^{*}_h}}{\varphi_h},$$\\ for every $t\in [0,T]$ converges weakly (up to a subsequence),  in the space $ L^q(0,T;U)$ to a control u such that the solution of :
$$\overset{.}{y}=Ay+Bu, \ y(0)=y_0,$$\\
satisfies $y(T)=0$. For every $h\in (0,h_1)$, let $y_h(.)$ denote the solution of
$$\overset{.}{y_h}=A_hy_h+B_hu_h, \ y_h(0)=P_h y_0.$$\\
\ Then,
\begin{itemize}
\item $y_h(T)=-h^{\beta}\left\|{\varphi_h}\right\|^{p-2}\varphi_h$;
\item The sequence $(\tilde{P_h}{y_h})_{0<h<h_1}$ converges strongly (up to subsequence) in the space $L^q(0,T;X)$, to y(.).
\end{itemize}

\par Futhermore, there exists $M>0$ such that
\begin{equation}\int_0^T{\left\|{u(t)}\right|^p_U}\le M^{p/{(p-1)}} \left\|{y_0}\right|^{p/{(p-1)}}_X,\nonumber
\end{equation}
\ and, for every $h\in (0,h_1)$,
\begin{eqnarray}&& \int_0^T{\left\|{u_h(t)}\right\|^p_{U_h}}\le M^{p/{(p-1)}}\left\|{y_0}\right\|^{p/{(p-1)}}_X,\nonumber\\
&&h^{\beta}\left\|{\varphi_h}\right\|^p_{X_h}\le M^{p/{(p-1)}}\left\|{y_0}\right\|^{p/{(p-1)}}_X,\nonumber\\
&&\left\|{y_h(T)}\right\|_{X_h}\le M^{1/{(p-1)}}h^{\beta/p}\left\|{y_0}\right\|^{1/{(p-1)}}_X.
\end{eqnarray}
\end{Theorem}
\par \textbf{Remark 4.4} The left hand side of (47) is uniform observability type inequality for control system (2). This inequality is weaker than the uniform exact null controllability. No attempt has been made here to prove uniform exact null control for the approximation systems (2).
\par \textbf{Remark 4.5} A similar result holds if the control system (1) is exactly controllable in time T. However, due to assumption ($H_1$ ), the semigroup S(t) enjoys in general regularity properties. Therefore, the solution y(.) of the control system may belong to a subspace of X, whatever the control u is. For instance, in the case of the heat equation with a Dirichlet or Neumann boundary control, the solution is the smooth function of the state variable x, as soon as $t>0$ ,for every control and initial condition $y_0\in L^2$. Hence, exact controllability does not hold in this case $L^2$.

\par Moreover, one may wonder under which assumptions the control u is the control, is defined by (27), such that y(T)=0. As in \cite{LT06}, the following proposition give an answer:
\begin{Proposition} With the notations of theorem, if the sequence of real numbers $\left\|{\psi
_h}\right\|_{X_h}$, $0<h<h_1$, is moreover bounded, then the control u is the unique control, is defined by (27), such that y(T)=0. Moreover, the sequence $(\tilde{Q}_hu_h)_{0<h<h_1}$ converges strongly (up to a sequence) in the space $L^q(0,T;U) $ to the control u.
\par A sufficient condition on $y_0\in X$, ensuring the boundedness of the sequence $(\left\|{\varphi_h}\right\|_{X_h})_{0<h<h_1}$, is the following : there exists $\eta>0$ such that the control system $\overset{.}{y}=Ay+Bu$ is exactly null controllable in time t, for every $t\in [T-\eta, T+\eta]$, and the trajectory $t\mapsto S(t)y_0$ in X, for $t\in [T-\eta, T+\eta]$, is not contained in a hyperplane of X.
\par Other sufficient condition on control u, also ensuring the boundedness of the sequence $(\left\|{\varphi_h}\right\|_{X_h})_{0<h<h_1}$, is the following : there exists $\eta>0$ such that the control system $\overset{.}{y}=Ay+Bu$ is exactly null controllable in time t, for every $t\in [T-\eta, T+\eta]$, and with the control u is defined as (27), the trajectory $t\mapsto S(t-\xi)Bu(\xi)$ in X, for $t\in [T-\eta, T+\eta]$, every $\xi\in (0,t)$ is not contained in a hyperplane of X.
\end{Proposition}

\section{Proof of the main results}
\begin{enumerate}
\item \textbf{The proof of theorem}:
\par \begin {proof} 
\par For convenience, we first state the following useful approximation lemma, whose proof readily follows that of \cite{LT06}, \cite{LT00}. The proof of this lemma is provided in the Appendix.

\begin{Lemma}There exists $C_{9}>0$ such that, for all $t\in (0,T] $ and $h\in (0,h_0)$, there holds

\begin{equation}{\left\|{(e^{tA^{*}_h}{P_h}-P_h S(t)^{*})\psi}\right\|}_{X_h}\le {C_9}\frac{h^s}{t}{\left\|{\psi}\right\|}_X,\end{equation}

\begin{equation}\left\| {\tilde{Q_h}{B^{*}_h}{e^{tA^{*}_h}}{\psi_h}}\right\|_U\le \frac{C_9}{t^{\gamma}}{\left\|{\psi_h}\right\|}_{X_h},\end{equation}
\\ for every $\theta \in [0,1]$ .
\begin{equation}\left\|{\tilde{Q_h}{ {B^{*}_h}{e^{tA^{*}_h}}{\psi_h}}-{B^{*}S(t)^{*}\tilde{P_h}{\psi_h}}}\right\|_U\le C_9\frac{h^{s(1-\gamma)\theta}}{t^{\theta+(1-\theta)\gamma}}{\left\|{\psi_h}\right\|_{X_h}}~\textrm{every}~{\psi_h}\in {X_h}.\end{equation}
 
\end{Lemma}
\par  We carry out proving the theorem as follows:

\par The degree of unboundedness $\gamma$ of the control operator B is lower than $\frac{1}{p}$, there exists $\theta \in (0,1)$ such that $0<\theta +(1-\theta)\gamma<\frac{1}{p}$.
\par For all $h\in (0,h_0)$ and $\psi_h\in {X_h}$ we have
\begin{eqnarray}\int_0^T{\left\|{ \tilde{Q_h}{B^{*}_h}{e^{tA^{*}_h}}{\psi_h}}\right\|}^p_U dt&=&\int_0^T(\left\| {\tilde{Q_h}{B^{*}_h}{e^{tA^{*}_h}}{\psi_h}}\right\|^p_U-{\left\|{{B^{*}}{S(t)^{*}}{\tilde{P_h}}{\psi_h}}\right\|}^p_U)dt \nonumber\\
&&+\int_0^T{\left\|{{B^{*}}{S(t)^{*}}{\tilde{P_h}}{\psi_h}}\right\|}^p_U dt.\end{eqnarray}
\par We estimate two terms of right hand side of (53).

\par  The control system is exactly null controllable in time T, then there exists a positive real number $C>0$ such that
\begin{eqnarray}  \int_0^T  {\left\|{B^*S(t)^*\tilde{P_h}{\psi_h}}\right\|^p}dt \ge C \left\|{S(T)^*\tilde{P_h}{\psi_h}}\right\|^p_X.\end{eqnarray}

\par We have the following inequality
 \begin{equation}\left|{y^p-z^p}\right|<p(y^{p-1}+z^{p-1})\left|{y-z}\right|,\end{equation} where $y,z \in R^+, p>1$.
\par Indeed, we apply mean-value theorem for $f(x)=x^p (p>1, x\in R^+)$ , there exists $\xi\in (y,z)$ such that
\begin{eqnarray} \left|{y^p-z^p}\right|&=&\left|{f'(\xi)}\right|\left|{y-z}\right| \ \nonumber\\
&=&p\left|{\xi^{(p-1)}}\right|.\left|y-z\right|\nonumber\\
&<&p(y^{p-1}+z^{p-1}).\left|y-z\right| .\nonumber\end{eqnarray}

\par  We apply the above inequality and make use of (40), (28), (44), (50) to obtain
\begin{eqnarray} &&\left|{\left\|{P_h S(T)^*\tilde{P_h}\psi_h}\right\|^p_{X_h}-\left\|{e^{TA^*_h}{\psi_h}}\right\|^p_{X_h}}\right|\nonumber\\
&\le& p(\left\|{P_hS(T)^*\tilde{P_h}\psi_h}\right\|^{p-1}_{X_h}+\left\|{e^{TA^*_h}\psi_h}\right\|^{p-1}_{X_h})\nonumber\\
&& \ \ \ \ \ \ \ \ \ \ \ \ \ \ \ \ \ \ \ \ \ \ \ \ \ \ \ \ \ \ \ \times \left|{\left\|{P_hS(T)^*\tilde{P_h}\psi_h}\right\|_{X_h}-\left\|{e^{TA^*_h}\psi_h}\right\|_{X_h}}\right|\nonumber\\
&\le& p(C_5C_1e^{\omega t}\left\|{\psi_h}\right\|^{p-1}_{X_h}+C_7\left\|{\psi_h}\right\|^{p-1}_{X_h}).\left\|{P_hS(T)^*\tilde{P_h}\psi_h-e^{TA^*_h}\psi_h}\right\|_{X_h}\nonumber\\
&\le& C_p\left\|{\psi_h}\right\|^{p-1}_{X_h}C_9C_5h^s \left\|{\psi_h}\right\|_{X_h}\nonumber\\
&\le& C_{14} h^s \left\|{\psi_h}\right\|^p_{X_h}.\nonumber \end{eqnarray}

\par Therefore , from above estimate and (39), we get
\begin{equation}{\left\|{{e^{TA^{*}_h}}{\psi_h}}\right\|^p_{X_h}}-{C_{14}}h^s{{\left\|{\psi_h}\right\|}^p_{X_h}}\le {\left\|{P_hS(T)^{*}{\tilde{P_h}}{\psi_h}}\right\|}^p_{X_h}\le C_5^p{\left\|{S(T)^{*}{\tilde{P_h}}{\psi_h}}\right\|}^p_{X}.\end{equation}

\par Combine (54) with (56) we have:
\begin{equation}\int_0^T{{\left\|{{B^{*}}{S(t)^{*}}{\tilde{P_h}}{\psi_h}}\right\|}^{p}_U}dt\ge {C_{15}}{\left\|{{e^{TA^{*}_h}}{\psi_h}}\right\|^p_{X_h}}-{C_{14}}h^s{{\left\|{\psi_h}\right\|}^p_{X_h}}.\end{equation}

\par For the first term on the right hand side of (53), one has, using (33), (51), (52) and applying the inequality (55)

\begin{eqnarray}&&\left|{ {{\left\|{\tilde{Q_h} {B^{*}_h}{e^{tA^{*}_h}}{\psi_h}}\right\|}^p_U-{\left\|{B^{*}{S(t)^{*}}{\tilde{P_h}}{\psi_h}}\right\|}^p_U}}\right| \nonumber\\
&\le&p(\left\|{\tilde{Q_h} {B^{*}_h}{e^{tA^{*}_h}}{\psi_h}}\right\|^{p-1}_U+\left\|{B^{*}{S(t)^{*}}{\tilde{P_h}}{\psi_h}}\right\|^{p-1}_U)\nonumber\\
&& \ \ \ \ \ \ \ \ \ \ \ \ \ \ \ \ \ \ \ \ \ \ \ \ \ \ \ \times\left|{\left\|{\tilde{Q_h} {B^{*}_h}{e^{tA^{*}_h}}{\psi_h}}\right\|_U-\left\|{B^{*}{S(t)^{*}}{\tilde{P_h}}{\psi_h}}\right\|_U}\right|\nonumber\\
&\le&p(\frac{C^{p-1}_9}{t^{\gamma(p-1)}}\left\|{\psi_h}\right|^{p-1}_{X_h}+C^{p-1}_3\frac{e^{\omega t(p-1)}}{t^{\gamma(p-1)}}\left\|{\psi_h}\right\|^{p-1}_{X_h})\nonumber\\
&& \ \ \ \ \ \ \ \ \ \ \ \ \ \ \ \ \ \ \ \ \ \ \ \ \ \ \ \ \ \ \ \times \left|{\left\|{{ \tilde{Q_h}{B^{*}_h}{e^{tA^{*}_h}}{\psi_h}}-{{B^{*}}{S(t)^{*}}{\tilde{P_h}}{\psi_h}}}\right\|_U}\right|\nonumber\\ 
&\le&\frac{C_{16}}{t^{\gamma(p-1)}}\left\|{\psi_h}\right\|^{p-1}_{X_h}.C_9\frac{h^{s(1-\gamma)\theta}}{t^{\theta+(1-\theta)\gamma}}\left\|{\psi_h}\right\|_{X_h}\nonumber\\
&\le& C_{17} \frac{h^{s(1-\gamma)\theta}}{t^{\theta+(1-\theta)\gamma+\gamma(p-1)}}\left\|{\psi_h}\right\|^p_{X_h}.\nonumber \end{eqnarray}

\par We have $\gamma<\frac{1}{p}(p\ge 1)$, therefore $\theta+(1-\theta)\gamma+\gamma(p-1)<1$ and we can get, by integration,
\par $\left|{\int_0^T{({\left\|{ \tilde{Q_h}{B^{*}_h}{e^{tA^{*}_h}}
{\tilde{P_h}}{\psi_h}}\right\|}^p_U-{\left\|{{B^{*}}{S(t)^{*}}{\psi_h}}\right\|}^p_U)}dt}\right| \le { C_{18}h^{s(1-\gamma)\theta}}{\left\|{\psi_h}\right\|^{p}_{X_h}}.$

\par Therefore,
\begin{equation}\int_0^T{{\left\|{ \tilde{Q_h}{B^{*}_h}{e^{tA^{*}_h}}{\psi_h}}\right\|}^p_Udt} \ge \int_0^T{\left\|{{B^{*}}{S(t)^{*}}{\psi_h}}\right\|}^p_Udt-{ C_{18}h^{s(1-\gamma)\theta}}{\left\|{\psi_h}\right\|^{p}_{X_h}}.\end{equation}
\par We choose a real number $\beta $ such that $0\le {\beta}\le {s(1-\gamma)\theta}$. Combine results (53), (57), (58) we have inequality (47).\\

\par For $h\in (0,h_1)$, the functional $J_h$ is convex, and inequality (47), is coercive. Therefore, it admits a solution minimum at $\varphi_h\in {X_h}$ so that
\begin{equation} 0=\bigtriangledown{{J_h}(\varphi_h)}={G_h}(T){\varphi_h}+{h^{\beta}}{\left\|{\varphi_h}\right\|}^{p-2}{\varphi_h}+e^{TA_h}{P_h}{y_0},\nonumber\end{equation}
\par where ${G_h}(T)=\int_0^T{ {\left\|{B^{*}_h e^{tA^{*}_h}{\varphi_h}}\right\|}^{p-2}{e^{tA_h}B_hB^{*}_he^{tA^{*}_h}}dt}$ is the Gramian of the semidiscrete system.
\par With $u_h(t)=\left\|{B^{*}_h e^{(T-t)A^{*}_h}{\varphi_h}}\right\|^{p-2}{B^{*}_he^{(T-t)A^{*}_h}}{\varphi_h}$ is chosen then, the solution $y_h(.)$ satisfies

\begin{eqnarray}y_h(T)&=&e^{TA_h}{y_h(0)}+\int_0^T{e^{(T-t)A_h}B_hu_h(t)dt}\nonumber\\
&=&e^{TA_h}{P_h}{y_0}+{G_h}(T){\varphi_h}\nonumber\\
&=&-h^{\beta}{\left\|{\varphi_h}\right\|}^{p-2}{\varphi_h}.\nonumber\end{eqnarray}

\par Note that, since $J_h(0)=0$, there must hold, at the minimum, $J_h(\varphi_h)\le 0$. Hence, using the observability inequality (47) and the Cauchy-Schwarz inequality, one gets
\begin{eqnarray} c\left\|{e^{TA^*_h}\varphi_h}\right\|^p_{X_h}&\le& \int_0^T{\left\|{B^*_he^{tA^*_h}\varphi_h}\right\|^p_{U_h}}+h^{\beta}\left\|{\varphi_h}\right\|^p_{X_h}\nonumber\\
&\le& 2\left\|{e^{TA^*_h}\varphi_h}\right\|_{X_h}\left\|{P_h y_0}\right\|_{X_h}, \nonumber
\end{eqnarray}
 and thus,
\begin{equation} \left\|{e^{TA^*_h}\varphi_h}\right\|_{X_h}\le (\frac{2}{c})^{1/{(p-1)}}(\left\|{P_h y_0}\right\|_{X_h})^{1/{(p-1)}}.
\end{equation}
\par As a consequence,
\begin{equation}\int_0^T{\left\|{B^*_he^{tA^*_h}\varphi_h}\right\|^p_{U_h}}\le (\frac{2^p}{c})^{1/{(p-1)}}(\left\|{P_hy_0}\right\|_{X_h}^{p/{(p-1)}}),
\end{equation}
and $h^{\beta}\left\|{\varphi_h}\right\|^p_{X_h}\le (\frac{2^p}{c})^{1/{(p-1)}}(\left\|{P_hy_0}\right\|_{X_h}^{p/{(p-1)}})$, and the estimates (49) follow.
\end {proof} 
\item \textbf{Proof of proposition}
\begin {proof} 
\par If the sequence $(\left\|{\tilde{P}_h\varphi_h}\right\|_X)_{0<h<h_1}$ is bounded then up to a subsequence, it converges weakly to an element $\varphi\in X$. It follows from the estimate (52) that $u(t)=\left\|{B^*S(T-t)^*\varphi}\right\|^{p-2}B^*S(T-t)^*\varphi$ for every $t\in [0,T]$. Moreover, $\tilde{Q}_hu_h$ tends strongly to u in $L^q(0,T;U)$. Indeed, for $t\in[0,T]$,
\begin{eqnarray} &&\tilde{Q}_h u_h(t)-u(t)\nonumber\\
&=& \tilde{Q}_h \left\|{B^{*}_h e^{(T-t)A^{*}_h}{\varphi_h}}\right\|^{p-2}{B^{*}_he^{(T-t)A^{*}_h}}{\varphi_h}\nonumber\\
&& \ \ \ \ \ \ \ \ \ \ \ \ \ \ \ \ \ \ \ \ \ \ \ \ \ \ \ \ \ \ \ \ \ \ \ \ \ \ \ -\left\|{B^*S(T-t)^*\varphi}\right\|^{p-2}B^*S(T-t)^*\varphi\nonumber\\
&=&\left\|{B^{*}_h e^{(T-t)A^{*}_h}{\varphi_h}}\right\|^{p-2}( \tilde{Q}_h B^{*}_he^{(T-t)A^{*}_h}-B^*S(T-t)^*\tilde{P}_h)\varphi_h\nonumber\\
&&+\left\|{B^{*}_h e^{(T-t)A^{*}_h}{\varphi_h}}\right\|^{p-2} B^*S(T-t)^*(\tilde{P}_h\varphi_h-\varphi)\nonumber\\
&&+B^*S(T-t)^*\varphi(\left\|{B^{*}_h e^{(T-t)A^{*}_h}{\varphi_h}}\right\|^{p-2}-\left\|{B^*S(T-t)^*\varphi}\right\|^{p-2}) .\end{eqnarray}
\par Since the $\varphi_h$ are bounded, then the $\left\|{u_h}\right\|$ are bounded. From that, we imply the $\left\|{B^{*}_h e^{(T-t)A^{*}_h}{\varphi_h}}\right\|^{p-2}$ are bounded .
\par Using (52), the first term of right hand side of (61) tends to zero clearly. For the second term, for every $t\in[0,T]$ the operator $B^*S(T-t)^*$ is compact, as a strongly limit of finite rank operators and since   $\tilde{P}_h\varphi_h - \varphi$ tends to weakly to zero, it follows the second term of the right hand side (61) tends to zero. Furthermore, by applying inequality (55) we get
\begin{eqnarray}&& \left\|{B^{*}_h e^{(T-t)A^{*}_h}{\varphi_h}}\right\|^{p-2}-\left\|{B^*S(T-t)^*\varphi}\right\|^{p-2}\nonumber\\
&<&(2-p)(\left\|{B^{*}_h e^{(T-t)A^{*}_h}{\varphi_h}}\right\|^{p-3}+\left\|{B^*S(T-t)^*\varphi}\right\|^{p-3})\nonumber\\
&& \ \ \ \ \ \ \ \ \ \times (\left\|{B^{*}_h e^{(T-t)A^{*}_h}{\varphi_h}-B^*S(T-t)^*\varphi}\right\|).\nonumber\end{eqnarray}
\par As the $\left\|{B^{*}_h e^{(T-t)A^{*}_h}{\varphi_h}}\right\|^{p-3}$ are bounded and inequality (52) is used again. Hence, the third term tends to zero clearly.
\par The control u is such that y(T)=0, hence the vector $\varphi$ must be solution of $\nabla J^*(\psi)=0$, where J is defined as in Theorem 2. Since J is convex, $\varphi$ is the minimum of $J^*$, that is, u is the control such that $y(T)=0$. \\

\par We next prove, by contradiction, that the sufficient conditions provided in the statement of the proposition implies that the sequence $(\left\|{\varphi_h}\right\|_{X_h})_{0<h<h_1}$ is bounded. As the proof of the first sufficient condition is found in \cite{LT06}, we give here the proof only for the second sufficient condition. If the sequence $(\left\|{\varphi_h}\right\|_{X_h})_{0<h<h_1}$ is not bounded, then, up to subsequence, $\tilde{P_h}(\varphi_h/ \left\|{\varphi_h}\right\|_{X_h})$ converges weakly to $\Phi$ in X, as h tends to 0. For every $t\in [T-\eta,T+\eta]$, the control system is exactly null controllable in time t; and thus, from (60), the sequence $\int_0^t{<B^*_h e^{(t-s)A^*_h}\varphi_h,Q_h u(s)>_{U_h}ds}$ is bounded, uniformly for $h\in(0,h_1)$. Thus, passing to the limit, one gets 
\begin{equation}\int_0^t{<\Phi, S(t-s)Bu(s)>_X ds}=0.\nonumber\end{equation}
\par  This equality is equivalent to the fact that : there exists $\xi\in(0,t)$ such that  $<\Phi, S(t-\xi)Bu(\xi)>_X=0$. This contradicts the fact that the trajectory $t\mapsto S(t-\xi)Bu(\xi)$, $t\in[T-\eta,T+\eta]$ and every $\xi\in(0,t)$, is not contained in a hyperplane of X.
\end {proof} 
\end{enumerate}

\section{Numerical simulation for the heat equation with Dirichlet boundary control}
\par In this section, we give an example of a situation where the theorem 3.1 are satisfied.
\par Let $\Omega \subset R^n$ be an open bounded domain with sufficiently smooth boundary $\Gamma$. We consider the Dirichlet mixed problem for the heat equation:
\begin{equation}\overset{.}{y}= \Delta y+c^2 y  ~\textrm{in}~ (0,T) \times \Omega \nonumber\end{equation}
\begin{equation} y(0,.)=y_0   ~\textrm{in}~ \Omega\nonumber\end{equation}
\begin{equation} y=u          ~\textrm{in}~  (0,T)\times \Gamma=\Sigma,\nonumber\end{equation}
\\with boundary control $u\in L^6(0,T;L^2(\Gamma))$ and $y_0\in L^2(\Omega)$.
\par Set $X=L^2(\Omega)$ and $U=L^2(\Gamma)$. We introduce the self-adjoint operator:
\begin{equation} Ah=\Delta h+c^2 h:D(A)=H^2\cap H^1_0 \to L^2(\Omega) .\nonumber\end{equation}

\par The adjoint $B^*\in L(D(A^*),U)$ of B is given by
\begin{equation} B^* \psi=-\frac{\partial \psi}{\partial \nu}, \psi \in D(A^*).\nonumber\end{equation}
\par Moreover, the degree of unboundedness of B is $\gamma =\frac{3}{4}+\epsilon$ ($\epsilon>0$) (see \cite{LT00}, section 3.1). 
\begin{enumerate}

\item \textbf{One-dimensional Finite-Difference semi-discretized model:}
\par We next introduce a semi-discretized model of the above heat equation, using 1D Finite-Difference.
\par For simplicity, we set $\Omega=(0,1)$, $\Gamma=\left\{0,1\right\}$, c=1 and T=1.
\par Given $n\in \mathbb{N}$ we define $h=\frac{1}{n+1}>0$. We consider the following simplex mesh:
\begin{center}
$\Omega_h=\left\{{x_0=0; x_i=ih, i=1,...,N; x_{n+1}=1}\right\},$ 
\end{center}
\par which divides [0,1] into n+1 subintervals $I_j=[x_j,x_{j+1}]$ j=0,...,n+1. Set
\begin{center}
 $ X_h=\left\{{y\in C^0(\Omega_h)}\right\},$
\end{center}
\begin{center}
$U_h=\left\{{y\in C^0(\Gamma)}\right\}$.
\end{center}
\par Define $\tilde{P}_h$ (resp., $\tilde{Q}_h$) as the canonical embedding from $X_h$ into $D((-A)^{1/2})$ (resp., from $U_h $ to U). For $x_h\in X_h$ and $u_h\in U_h$, set, $\tilde{P}_h(x_h)=x_h$ and $\tilde{Q}_h(u_h)=u_h$. For $y\in D((-A)^{1/2})'=H^1(\Omega)'$, set $P_h y=(y_1,..,y_i,..,y_{n+2})$ where $y_i=y((i-1)h)$ and, for $u\in U$, set, $Q_h u=(u_1,..,u_i,..,u_{n+2})$ where $u_i=u((i-1)h)$.
\par It is clear that the assumptions $(H_{3.1})$ and $(H_{3.3})$ are here satisfied. Our aim is next  to verify the inequalities in $(H_{3.2})$ and $(H_{3.4})$.\\
\par In order to get these inequalities, it will necessary to making use of the following usual approximation properties (see \cite{LT00}, section 5):
\par (i) $\left\|{\Pi_h y-y}\right\|_{H^l(\Omega)}\le ch^{s-l}\left\|{y}\right\|_{H^s(\Omega)}$, $s\le r+1, \  \ s-l\ge 0,\ \ \  0\le l\le 1$,
\\and the inverse approximation properties
\par (ii) $\left\|{y_h}\right\|_{H^{\alpha}(\Omega)}\le ch^{-\alpha}\left\|{y_h}\right\|_L^2(\Omega)$, $0\le \alpha\le 1$.
\par (iii) $h^{-1}\left\|{y-\Pi_h y}\right\|_{L^2_{(\Gamma)}}+\left\|{(I-\Pi_h  )\frac{\partial y }{\partial \nu}}\right\|_{L^2(\Gamma)}\le c h^{s-\frac{3}{2}}\left\|{y}\right\|_{H^s(\Omega)}$, $\frac{3}{2}<s<r+1, y\in H^s(\Omega)$.
\par (iv) $\left\|{y_h}\right\|_{L^2(\Gamma)}+h \left\|{\frac{\partial y_h}{\partial \nu}}\right\|_{L^2(\Gamma)}\le Ch^{-\frac{1}{2}}\left\|{y_h}\right\|_{L^2(\Omega)}$, $y_h\in V_h$.
\\where r is the order of approximation (degree of polynomials) and $\Pi_h$ is the orthogonal projection of $L^2(\Omega)$ onto $V_h$ .\\
\par First, by applying (i) we easily get the inequality (36)
\begin{eqnarray} \left\|{(I-\hat{A^*}^{-\gamma+\frac{1}{2}}\tilde{P}_h P_h)\psi}\right\|_{L^2(\Omega)}&\le& C h^2\left\|{\psi}\right\|_{H^2(\Omega)}\nonumber\\
&\le&C h^2\left\|{\psi}\right\|_{D(A^*)}\nonumber\\
&\le&C h^2\left\|{A^*\psi}\right\|_X,\nonumber\end{eqnarray}
\par in this case s=2.\\
\par We next verify the inequality (37) as follows
\begin{eqnarray}&&\left\|{((-\hat{A}^{*})^{\gamma})}{\left({I-{{\hat{(A^*)}^{-\gamma+\frac{1}{2}}}\tilde{P_h}{P_h}}}\right)\psi}\right\|_X \nonumber\\
&\le & C \left\|{ ( I-{\hat{(A^*)}^{-\gamma+\frac{1}{2}}\tilde{P_h}{P_h}})\psi }\right\|_{D((-\hat{A}^{*})^{\gamma})} \nonumber\\
&\le & Ch^{s-l}\left\|{\psi}\right\|_{D(A^*)}\nonumber\\
&\le &C{h^{s({1-\gamma})}}\left\|{A^{*}\psi}\right\|_X, \nonumber \end{eqnarray}
\\where we have used (i) with s=2, $D(A^*)=H^s(\Omega)$ and $D((-\hat{A}^{*})^{\gamma})=H^l(\Omega)$ .\\
\par For the inequality (39), we apply (iii) with s=2 as 
\begin{eqnarray} \left\|{(I-\tilde{Q}_h Q_h)B^*\psi}\right\|_{L^2(\Gamma)}& =&\left\|{(I-\tilde{Q}_h Q_h)\frac{\partial \psi}{\partial \nu}}\right\|_{L^2(\Gamma)}\nonumber\\
&\le&Ch^{1/2}\left\|{\psi}\right\|_{H^2(\Omega)}\nonumber\\
&\le&Ch^{s(1-\gamma)}\left\|{\psi}\right\|_{D(A^*)}\nonumber\\
&\le&Ch^{s(1-\gamma)}\left\|{A^*\psi}\right\|_X.\nonumber
\end{eqnarray}\\
\par For the inequality (43), by using (iv) and (40) we get
\begin{eqnarray}\left\|{B^{*}{\hat{(A^*)}^{-\gamma+\frac{1}{2}}}\tilde{P_h}\psi_h}\right\|_U&=&\left\|{\frac{\partial( {\hat{(A^*)}^{-\gamma+\frac{1}{2}}}\tilde{P_h}\psi_h)}{\partial \nu}}\right\|_{L^2(\Gamma)}\nonumber\\
&\le&C h^{-\frac{3}{2}}\left\|{{\hat{(A^*)}^{-\gamma+\frac{1}{2}}}\tilde{P_h}\psi_h}\right\|_{L^2(\Omega)}\nonumber\\
&\le& Ch^{-\frac{3}{2}}\left\|{\psi_h}\right\|_{X_h}.\nonumber\end{eqnarray}
\par Therefore, the inequality (43) is satisfied for s=2, $\gamma=\frac{3}{4}+\epsilon$.\\
\par  Moreover, the assumption $(H_{4.2})$ is satisfied with s=2 (see \cite{LT00}).
\par  Hence, theorem 3.1 applies, with $\beta=0.16 $, for instance.\\
\par We then consider the following finite difference approximation of the above heat equation as follws
\begin{center}
\begin{eqnarray}&& \dot{y}_j=\frac{1}{h^2}[y_{j+1}+y_{j-1}-2y_j]+c^2y_j \,  0<t<T,  j=1,...,n\nonumber\\
&&y_j(0)=y_{j0}                                                               \ ,  j=1,...,n\nonumber\\
&&y_0(t)=y_{n+1}(t)= u_h                                     \  ,    0<t<T,\nonumber \end{eqnarray}
\end{center}
where $y\in R^{n+2}$ , $y_0\in R^{n+2}$ , $u_h\in R$ and
\begin{displaymath}
\mathbf{A_h} =\frac{1}{h^2}
\left( \begin{array}{cccccc}
0 & 0 &0& \ldots&0 \\
0 & (c^2h^2-2)&1 & \ldots&0 \\
0&1&(c^2h^2-2)&\ldots&0\\
%0&0&1&(c^2h^2-2)&\ldots&0\\
\vdots & \vdots & \ddots&\vdots&\vdots\\
%0&\ldots&(c^2h^2-2)&1&0&0\\
0&\ldots&(c^2h^2-2)&1&0\\
0&\ldots&1&(c^2h^2-2)&0\\
0 & \ldots&0 &0&0
\end{array} \right),
\end{displaymath}

\begin{displaymath}
\mathbf{B_h} =
\left( \begin{array}{ccc}
1\\
0\\
\vdots \\
0\\
1
\end{array} \right).
\end{displaymath}

\item\textbf{ Numerical simulation}
\par The minimization procedure described in Theorem 3.1 has been implement for d=1, by using a simple gradient method that has the advandtages not to require the complex computations and this method can applied with any power p. However, the computation of gradient of $J_h$ is very expensive since the gradient is related to the Gramian matrix.

\begin{center}
\begin{tabular}{|l|l|l|l|}

\hline name&$S_h$ & h &$y_0$
\\  \hline \hline
1D-10&10&$10^{-1}$&$y_0(x)=e^{-x^2}$ \\
1D-100&100&$10^{-2}$&$y_0(x)=e^{-x^2}$ \\
1D-500&500& $2.10^{-3}$& $y_0(x)=e^{-x^2}$ \\ \hline
\end{tabular}
\par Table 1: Data for the one-dimensional heat equation.
\end{center}

\begin{center}
\begin{tabular}{|l|l|l|l|l|}

\hline name&$\left\|{\varphi_h}\right\|_X$ & $h^{\beta} (\beta=0.16)$ &$\left\|{y_h(T)}\right\|$
\\  \hline \hline
1D-10&0.1690&$0.6814$&0.4775 \\
1D-100&0.7960&$0.4779$&0.4565 \\
1D-500&2.0570&$0.3699$&0.4273\\ \hline
\end{tabular}
\par Table 2: Numerical results for one dimensional equation for $\beta=0.16$.
\end{center}

\begin{center}
\begin{tabular}{|l|l|l|l|l|}
\hline name&$\left\|{\varphi_h}\right\|_X$ & $h^{\beta} (\beta=2)$ &$\left\|{y_h(T)}\right\|$
\\  \hline \hline
1D-10&4.4266&$10^{-2}$&0.0111\\
1D-100&4.8933&$10^{-4}$&1.3467e-004 \\
1D-500&5.0956&$4.10^{-6}$&5.5178e-006\\ \hline
\end{tabular}
\par Table 3: Numerical results for one dimensional equation for $\beta=2$.
\end{center}

\par Numerical simulation are carried out with a space-discretization step equal to 0.005, with the data of Table 1. The numerical results are provided in Table 2 for beta =0.16 and Table 3 for beta=2.
\par The convergence of the method is very slow. From the result of Theorem 3, the final state $y_h(T)$ is equal to $-h^{\beta}\left\|{\varphi_h}\right\|^{p-2}\varphi_h$ in which $\varphi_h$ is minimizer of $J_h$. We note that the maximum value for which the theorem asserts the convergence is very small. For such a small value of $\beta$ (for instance $\beta=0.16$), $h^{\beta}$ converges to 0 very slowly. It follows that $y_h(T)$ converges very extremely slow.
\par In practice, the unique minimizer $\varphi_h$ of $J_h$ is obtained through the simple gradient method in which the step is equal to 0.01 and the error $\epsilon=10^{-2}$ is taken. With the small value $\beta=0.16$, it took a long time to achieve the results in Table 2. Namely, for n=500, it took more one month to get the result after 1000 iterations. It is clearly seen from Table 2 that the convergence of $y_h(T)$ is very slow. These results illustrate the difficult in using our method to compute controls. Although the case beta=2 is not covered by our main theorem, the method seems to converge for this value for beta and we provide hereafter the numerical results in Table 3. Since the value of the beta is greater, the convergence is quicker.
% In order to observe the convergence more easily, we will take $\beta=2$ in pratice. Namely, in Table 2, it is easily seen that we collect the result by the conjudate gradient method (with $\beta= 2$ is taken ). 

\end{enumerate}
\textbf{Appendix: proof of lemma}

\par \begin {proof}

\begin{itemize}

\item First of all, we will  prove (51)

\par For every $\psi \in D(A^*)$, one has\\
\begin{eqnarray}&&\left\|{\tilde{Q_h}B^*_h e^{tA^*_h}P_h{\psi}-B^*{\hat{(A^*)}^{-\gamma+\frac{1}{2}}}\tilde{P_h}{P_h}S(t)^*{\psi}}\right\|_U \le{\left\|{\tilde{Q_h}B^*_h e^{tA^*_h}P_h{\psi}}\right\|_U}\nonumber\\
&&\ \ \ \ \ \ \ \ \ \ \ \ \ \ \ \ \ \ \ \ \ \ \ \ \ \ \ \ \ \ \ \ \ \ \ \ \ \ \ \ \ \ \  \ \ \  \ \ \ \ \ \ \ \  +{\left\|{B^*{\hat{(A^*)}^{-\gamma+\frac{1}{2}}}\tilde{P_h}{P_h}S(t)^*{\psi}}\right\|_U}.\nonumber\\
\end{eqnarray}
\par We estimate each term of the right hand side of (62). Since ${B^{*}_h}={Q_h}{B^{*}}{\hat{(A^*)}^{-\gamma+\frac{1}{2}}}{\tilde{P_h}}$  and thus, using (40)
(42) (44) one gets

\begin{eqnarray} {\left\|{\tilde{Q_h}B^*_h e^{tA^*_h}{P_h}{\psi}}\right\|_U}
 &\le& C_5{\left\|{B^*{\hat{(A^*)}^{-\gamma+\frac{1}{2}}}\tilde{P_h}e^{tA^*_h}{P_h}{\psi}}\right\|_U}\nonumber\\
 &\le& C_5 C_6 h^{-\gamma s}{\left\|{e^{tA^*_h}{P_h}{\psi}}\right\|_{X_h}}\nonumber\\
  &\le& C^2_5 C_6 C_7 h^{-\gamma s} e^{\omega t}{\left\|{\psi}\right\|_X} . \end{eqnarray}

\par On the other hand, from (28), (40), (42) one obtains

\begin{eqnarray}\left\|{B^*{\hat{(A^*)}^{-\gamma+\frac{1}{2}}}\tilde{P_h}{P_h}S(t)^*{\psi}}\right\|_U  &\leq& C_6 h^{-\gamma s}{\left\|{P_hS(t)^*\psi}\right\|_{X_h}} \nonumber \\
&\leq & C_5 C_6 h^{-\gamma s} {\left\|{S(t)^*\psi}\right\|_{X}}\nonumber \\
& \leq & C_1 C_5 C_6 h^{-\gamma s} e^{\omega t}{\left\|{\psi}\right\|_{X}}.
\end{eqnarray}

\par Hence, combine (63), (64) with (62), there exists $C_{10} >0$ such that
\begin{equation}\left\|{\tilde{Q_h}B^*_h e^{tA^*_h}P_h{\psi}-B^*{\hat{(A^*)}^{-\gamma+\frac{1}{2}}}\tilde{P_h}{P_h}S(t)^*{\psi}}\right\|_U \le C_{10} h^{-\gamma s}{\left\|{\psi}\right\|_{X}}.\end{equation}
 for every $\psi \in D(A^*)$, every $t\in [0,T]$, and every $h\in (0,h_0)$.
\par Moreover, we get another estimate of this term. By using (33), (38), (40) , (42), (50) one gets

\begin{eqnarray} && \left\|{\tilde{Q_h}B^*_h e^{tA^*_h}P_h{\psi}-B^*{\hat{(A^*)}^{-\gamma+\frac{1}{2}}}\tilde{P_h}{P_h}S(t)^*{\psi}}\right\|_U\nonumber\\
& =& {\left\|{\tilde{Q_h}{Q_h}B^* {\hat{(A^*)}^{-\gamma+\frac{1}{2}}}{\tilde{P_h}}e^{tA^*_h}P_h{\psi}-B^*{\hat{(A^*)}^{-\gamma+\frac{1}{2}}}\tilde{P_h}{P_h}S(t)^*{\psi}}\right\|_U}\nonumber\\
&\le& \left\|{\tilde{Q_h}{Q_h}B^*{\hat{(A^*)}^{-\gamma+\frac{1}{2}}}{\tilde{P_h}} (e^{tA^*_h}P_h{\psi}-{P_h}S(t)^*{\psi})}\right\|_U\nonumber\\
& &+\left\|{\tilde{Q_h}{Q_h}B^* ({{\hat{(A^*)}^{-\gamma+\frac{1}{2}}}\tilde{P_h}}{P_h}-I)S(t)^*{\psi}}\right\|_U\nonumber\\
& &+\left\|{(\tilde{Q_h}{Q_h}-I)B^*S(t)^*{\psi}}\right\|_U\nonumber\\
&&+\left\|{B^*(I-{\hat{(A^*)}^{-\gamma+\frac{1}{2}}}\tilde{P_h}P_h)S(t)^*{\psi}}\right\|_U \nonumber\\
 &\le& C_5 C_6 h^{\gamma s}{\left\|{e^{tA^*_h}P_h{\psi}-{P_h}S(t)^*{\psi}}\right\|_{X_h}}\nonumber\\
&&+ C_5 C_3 {\left\|{(-\hat{A})^{\gamma} ({{\hat{(A^*)}^{-\gamma+\frac{1}{2}}}\tilde{P_h}}{P_h}-I)S(t)^*{\psi}}\right\|_X}\nonumber
\\ &&+C_4 h^{s(1-\gamma)}{\left\|{A^*S(t)^*{\psi}}\right\|_X}\nonumber
\\ &&+C_3 {\left\|{(-\hat{A})^{\gamma} ({{\hat{(A^*)}^{-\gamma+\frac{1}{2}}}\tilde{P_h}}{P_h}-I)S(t)^*{\psi}}\right\|_X} \nonumber
\\& \le& C_5 C_6 C_{9} \frac{h^{s(1-\gamma)}}{t}{\left\|{\psi}\right\|_X}+(C_3(C_5+1)+1)C_4 h^{s(1-\gamma)}{\left\|{A^*S(t)^*{\psi}}\right\|_X}\nonumber
\\ &\le& C_{11}{\frac{h^{s(1-\gamma)}}{t}}{\left\|{\psi}\right\|_X} .\end{eqnarray}

\par Then, raising (65) to the power $1-\gamma $, (66) to power to $\gamma$ and multiplying both result estimates, we obtain 
 $$\left\|{\tilde{Q_h}B^*_h e^{tA^*_h}P_h{\psi}-B^*{\hat{(A^*)}^{-\gamma+\frac{1}{2}}}\tilde{P_h}{P_h}S(t)^*{\psi}}\right\|_U \le \frac{C_{12}}{t^\gamma}{\left\|{\psi}\right\|_{X}}.$$
\par Hence,
\begin{equation}\left\|{\tilde{Q_h}B^*_h e^{tA^*_h}P_h{\psi}}\right\|_U \le \frac{C_{12}}{t^\gamma}{\left\|{\psi}\right\|_{X}}+{\left\|{B^*{\hat{(A^*)}^{-\gamma+\frac{1}{2}}}\tilde{P_h}{P_h}S(t)^*{\psi}}\right\|_U}.\end{equation}
\par From (29), (33), (36) one yieds
\begin{eqnarray}\left\|{B^*{\hat{(A^*)}^{-\gamma+\frac{1}{2}}}\tilde{P_h}{P_h}S(t)^*{\psi}}\right\|_U
&\le&\left\|{B^*(I-{\hat{(A^*)}^{-\gamma+\frac{1}{2}}}\tilde{P_h}P_h)S(t)^*{\psi}}\right\|_U +{\left\|{B^*S(t)^*\psi}\right\|_U}\nonumber\\
&\le& C_{13} \frac{e^{\omega t}}{t^\gamma}{\left\|{\psi}\right\|_X}.\end{eqnarray}

\par Combine (67) with (68)  and by setting $\psi= \tilde{P_h}{\psi_h}$ we get (51).

\item Finally, we prove (52). On the one hand, reasoning as above for obtaining (66), we get
\begin{equation}{\left\|{\tilde{Q_h}B^*_h e^{tA^*_h}P_h{\psi}-B^*S(t)^*{\psi}}\right\|_U} \le C \frac{h^{s(1-\gamma)}}{t}{\left\|{\psi}\right\|_X},\end{equation}  
\ for every $\psi \in D(A^*)$, every $t\in [0,T]$ and every $h\in (0,h_0)$.
 \par On the other hand, from (51) and setting $\psi= \tilde{P_h}{\psi_h}$ one obtains
\begin{eqnarray}\left\|{\tilde{Q_h}B^*_h e^{tA^*_h}P_h{\psi}-B^*S(t)^*{\psi}}\right\|_U
&\le&{\left\|{\tilde{Q_h}B^*_h e^{tA^*_h}P_h{\psi}}\right|_U}+{\left\|{B^*S(t)^*{\psi}}\right\|_U}\nonumber\\
&\le& \frac{C_9}{t^\gamma}\left\|{\psi_h}\right\|+C_3{\left\|{(-\hat{A^*})^{\gamma}S(t)^*{\psi}}\right\|_X}\nonumber\\
&\le& \frac{C}{t^\gamma}\left\|{\psi_h}\right\|_X.\end{eqnarray}

\ Raising (69) to the power $\theta$, (70) to the power $1-\theta$ and multiplying both resulting estimates, we obtain (52).

\par The proof of the inequality (50) is found in \cite{LT06}, \cite{LT00}.

\end{itemize}

\end {proof}

\section{Conclusion}
%\begin{itemize}

%\item 
\par We have shown that the appropriate duality techniques can be applied to solve (3), namely the Fenchel-Rockafellar theorem.
\par Additionally, it is also stated that under standard assumptions on the discretization process, for an exactly null controllable linear control system, if the semigroup of approximating system is uniformly analytic, and if the degree of unboundedness of the control operator is greater than $\frac{1}{2}$ then the unform observability type inequality is proved. Consequently, a minimization procedure was provided to build the aproximation controls. This is implemented in the case of the one dimensional heat equation with Dirichlet boundary control.  %The authors in [13] dealt with in case the degree of boundedness of the control operator is lower than $\frac{1}{2}$ ..% From that, we show  the existence of minimum of $\int_0^T{\left\|{u(t)}\right\|^q}dt$ $(q>2)$. Moreover, the authors in [5], gave the result for the existence of minimum of this function with $1<q\le 2$. Therefore, we achieve the conclusion for the existence of minimum of $\int_0^T{\left\|{u(t)}\right\|^q}dt$ with every $q>1$.
\par Note that, we only stress our problem on the case $\gamma\ge 1/2$. Some relevant problems for which $\gamma<1/2$ that are refered to \cite{LT06}.
\par One open question is given: how the above results change if we remove the assumption of uniform analyticity of the discretized semigroup.
%\end{itemize} 

\section*{Acknowledgement}
The author warmly thank Emmanuel Trelat for his helps and remarks. And the author also wishes to acknowledge Region Centre for its financial support.
\renewcommand{\refname}{References}

% Ket thuc van ban
\end{document}